\documentclass[11pt]{article}
\usepackage{amssymb}

\usepackage{amsfonts}
\usepackage{amsmath}

\newtheorem{thm}{Theorem}

\newtheorem{cor}[thm]{Corollary}

\newtheorem{defn}[thm]{Definition}
\newtheorem{exmp}[thm]{Example}

\newtheorem{lemma}[thm]{Lemma}

\newtheorem{prop}[thm]{Proposition}
\newtheorem{rem}[thm]{Remark}

\newenvironment{pf}[1][Proof]{\noindent\textbf{#1.} }{\ \rule{0.5em}{0.5em}}

\def\one{1\!\!1}

\begin{document}

\title{$WI$-posets, graph complexes \\ and $\mathbb{Z}_2$-equivalences}
\author{Rade T. \v Zivaljevi\' c\\
Mathematics Institute SANU, Belgrade}

\date{May 3, 2004}

\maketitle
\begin{abstract} An evergreen theme in topological graph
theory is the study of graph complexes, \cite{BobsonKozlov2}
\cite{Lovasz} \cite{Matousek} \cite{MatZieg}. The majority of
these complexes are $\mathbb{Z}_2$-spaces and the associated
$\mathbb{Z}_2$-index $Ind_{{\mathbb{Z}_2}}(X)$ is an invariant of
great importance for estimating the chromatic numbers of graphs.
We introduce $WI$-posets (Definition~\ref{def_weakly}) as
intermediate objects and emphasize the importance of Bredon's
theorem (Theorem~\ref{thm_bredon}) which allows us to use standard
tools of topological combinatorics for comparison of
$\mathbb{Z}_{2}$-homotopy types of $\mathbb{Z}_2$-posets. Among
the consequences of general results are known and new results
about $\mathbb{Z}_{2}$-homotopy types of graph complexes. It turns
out that, in spite of great variety of approaches and definitions,
all graph complexes associated to $G$ can be viewed as avatars of
the same object, as long as their $\mathbb{Z}_2$-homotopy types
are concerned. Among the applications are a proof that each
finite, free $\mathbb{Z}_2$-complex is a graph complex and an
evaluation of $\mathbb{Z}_2$-homotopy types of complexes
$Ind(C_n)$ of independence sets in a cycle $C_n$.
\end{abstract}

\section*{Introduction\label{sec_introduction}}


By a deep observation of L. Lov\'{a}sz \cite{Lovasz}, the
chromatic number $\chi (G)$ of a graph can be approximated from
below by integers reflecting the topological complexity of
associated \emph{graph complexes.} The impact of this observation
can not be overestimated and this direction of topological graph
theory has been for decades a vital part of topological
combinatorics, see \cite{Matousek} \cite{MatZieg}\ and the
references therein.

The first in the series of \ graph complexes is the so called \emph{%
neighborhood complex }$N(G):=\left\{ S\subset V_{G}\mid CN(S)\neq
\varnothing \right\} $, where $CN(S)$ is the set of all common neighbors of $%
S$ in $G$. Currently there exist a dozen of graph complexes, see
the references \cite{AlonFranklLovasz}
\cite{BabsonKozlov}--\cite{BobsonKozlov2}
\cite{Csorba}--\cite{Csorba+3} \cite{Kriz}--\cite{MatZieg}
\cite{MilgZveng} \cite{Sarkaria} \cite{Walker}. Many of them
originated from the neighborhood complex $N(G)$ and all of them
are used to produce lower bounds for the chromatic number of $G$
in terms of other numerical invariants. A central among these
invariants is the equivariant index $Ind_{\mathbb{Z}_{2}}(K)$,
which applies to graph complexes $K$ with fixed point free
involutions $\omega :K\longrightarrow K$.
$Ind_{\mathbb{Z}_{2}}(K)$ is defined as the minimum integer $n$
such that there exists a $\mathbb{Z}_{2}$-equivariant map $f
:K\longrightarrow S^{n}$. This integer is an invariant of the
$\mathbb{Z}_{2}$-homotopy type of the $\mathbb{Z}_{2}$-complex $K$
and it is not a surprise that much of the current research is
focused on clarifying the mutual relationship of different graph
complexes \cite{Csorba2004} \cite{Csorba+3}  \cite{Matousek}
\cite{MatZieg}.

In this paper we develop a unified approach to the problem of
comparing $\mathbb{Z}_{2}$-homotopy types of graph complexes. The
so called $WI$-posets (Definition~\ref{def_weakly}) are designed
to capture the essential features of the neighborhood complex
(lattice) $N(G)$ and to serve as a basis for construction of
\emph{graph posets.} The idea to use posets (lattices) as
intermediate objects in the construction of graph complexes is not
new. J.W.\ Walker introduced \emph{ortholattices} in \cite{Walker}
precisely for this purpose and was the first to emphasize the
functoriality of such a construction. The novelty of our approach
is in the systematic use of Bredon's theorem
(Theorem~\ref{thm_bredon}) which allows us to shift from
$\mathbb{Z}_{2}$-homotopy types to the ordinary homotopy types of
posets. This change of perspective brings in powerful and elegant
tools of topological combinatorics, notably Quillen fiber theorem
and its relatives. It also accounts for the greater generality and
conceptual simplicity achieved by the introduction of $WI$-posets.
Using this approach we obtain new $\mathbb{Z}_{2}$-homotopy
equivalences between graph complexes and posets (Sections
\ref{sec_all_equivalent} and \ref{sec_avatars}) and as a
consequence derive new conceptual proofs of related results of
Matou\v{s}ek and Ziegler \cite{MatZieg}, Csorba et al.
\cite{Csorba} \cite{Csorba+3}, and Lov\'{a}sz (\cite{MatZieg}
Sect.\ 5). Among the highlights are a proof
(Section~\ref{sec_which}) of the fact that each finite, free
$\mathbb{Z}_2$-complex is a graph complex (earlier proved by
Csorba \cite{Csorba2004}) and an analysis of the
$\mathbb{Z}_2$-homotopy types of complexes $Ind(C_n)$, yielding
results originally used by Bobson and Kozlov \cite{BobsonKozlov2}
in their solution of Lov\'{a}sz conjecture
(Section~\ref{sec_ind}).

The notation used in the paper is standard \cite{Bjorner}.
$G=(V_{G},E_{G})$ is a finite graph with $V_{G}$ and $E_{G}$ as
the sets of vertices and edges. All graphs are simple and
undirected. The collection of all chains in
a (finite) posets $P$ forms a simplicial complex called the order complex $%
\Delta (P)$ of $P$. A $\mathbb{Z}_{2}$-space is a topological
space $X$ equipped with a continuous involution $\omega
:X\rightarrow X,\, \omega ^{2}=1_{X}$. A
$\mathbb{Z}_{2}$-equivariant map $f:X\longrightarrow Y$ between
two $\mathbb{Z}_{2}$-spaces $X$ and $Y$ is a continuous map
satisfying the condition $f(\omega x)=\omega f(x)$. A $\mathbb{Z}_{2}$%
-equivariant map, or a $\mathbb{Z}_{2}$-map for short, is a $\mathbb{Z}_{2}$%
-equivalence if there exists a $\mathbb{Z}_{2}$-map $g:Y\longrightarrow X$
such that $g\circ f$ is $\mathbb{Z}_{2}$-homotopic to $1_{X}$ and $f\circ g$
is $\mathbb{Z}_{2}$-homotopic to $1_{Y\text{.}}$ A general reference for $G$%
-spaces, $G$-equivariant maps and related concepts and facts is \cite%
{tomDieck}. Expositions oriented towards applications in combinatorics can
be found in \cite{Matousek}, \cite{elisat} and \cite{guides}.

\bigskip

\section{Involutive and weakly involutive
posets}\label{sec_involutive}

\bigskip

\begin{defn}
\label{def_involutive} A poset $(Q,\leq )$ is involutive ($I$-poset) if it
is equipped with an involution $C:Q\longrightarrow Q$ which is either
monotone or antitone, i.e. which satisfies either the condition $x\leq
y\Rightarrow C(x)\leq C(y)$ or the dual condition $x\leq y\Rightarrow
C(x)\geq C(y)$. We also say that $(Q,\leq )$ admits a $\mathbb{Z}_{2}$%
-action or that $(Q,\leq )$ is a $\mathbb{Z}_{2}$-poset.
\end{defn}

\begin{defn}
\label{def_weakly} A \emph{weakly involutive poset }$(P,C)$\emph{, }or a $WI$%
-poset for short, is a finite poset $P$ equipped with a function $%
C:P\rightarrow P$ such that%
\begin{equation}
x\leq y\Rightarrow C(y)\geq C(x)  \label{eqn_weak1}
\end{equation}%
\begin{equation}
x\leq C(C(x))=C^{2}x  \label{eqn_weak2}
\end{equation}
\end{defn}

\begin{rem}
The theories of antitone and monotone $I$-posets are similar but there are
also some important differences. For example\ only antitone $I$-posets are $%
WI$-posets in the sense of Definition \ref{def_weakly}. Some
results in the paper are sensitive to this difference so whenever
necessary, it will emphasized what kind of $I$-posets we are
dealing with.
\end{rem}

\begin{defn}
\label{def_Lov} Suppose that $(P,C)$ is a $WI$-poset. Then the
associated (antitone) involutive poset $\mathcal{L}(P)$ is a
subposet of $P$ defined by $\mathcal{L}(P)=\left\{ x\in P\mid
C^{2}x=x\right\} $.
\end{defn}

An easy consequence of equations (\ref{eqn_weak1}) and (\ref{eqn_weak2}) is
the equality $C^{3}x=Cx$ which implies that $\mathcal{L}(P)$ is non-empty
and that $C$, restricted to $\mathcal{L}(P)$, is a genuine antitone
involution turning $(\mathcal{L}(P),C)$ into an involutive poset in the
sense of Definition \ref{def_involutive}. The involutive poset $(\mathcal{L}%
(P),C)$ is often called the Lov\'{a}sz poset associated to $(P,C)$ for the
reasons explained in Section \ref{sec_applications}.

\begin{defn}
\label{def_box_poset} The \emph{box poset} $\mathfrak{B}(P)$
associated to a $WI$-poset $(P,C)$ is a subposet of $P\times P$
defined by
\begin{equation*}
\mathfrak{B}(P)=\mathfrak{B}(P,C)=\left\{ (x,y)\in P\times P\mid x\leq
Cy\And y\leq Cx\right\} .
\end{equation*}
\end{defn}
It is desirable to isolate the ``correct'' notion of a morphism of
$WI$-posets which would turn $P\mapsto \mathcal{L}(P)$ and
$P\mapsto \mathfrak{B}(P)$ into genuine functors. If $f :
(P,C)\rightarrow (Q,C)$ is a monotone map of $WI$-posets such that
$f(C(x))=C(f(x))$ then obviously there exists a monotone map
$\bar{f}: \mathcal{L}(P)\rightarrow \mathcal{L}(Q)$ of associated
Lov\'{a}sz posets. This condition is unfortunately too
restrictive. Here is a natural condition on a monotone map $f :
(P,C)\rightarrow (Q,C)$ of $WI$-posets guaranteeing that the
associated map $F: \mathfrak{B}(P)\rightarrow \mathfrak{B}(Q), \,
(x,y)\mapsto (f(x),f(y))$, is well defined and monotone.

\begin{defn}\label{def_wimorfizam}
A monotone map $f : (P,C)\rightarrow (Q,C)$ of $WI$-posets is a
$WI$-morphism if $f(C(x))\leq C(f(x))$ for each $x\in P$.
\end{defn}

\begin{defn}
\label{def_int_poset} The \emph{poset of intervals} ($Int(Q),\preccurlyeq )$%
, associated to a poset $(Q,\leq )$, is by definition
\begin{equation*}
Int(Q)=\left\{ (x,y)\in P\times P\mid x\leq y\right\}
\end{equation*}%
where $(x,y)\preccurlyeq (x^{\prime },y^{\prime })\Leftrightarrow
x\leq x^{\prime }\leq y^{\prime }\leq y.$ The elements of $Int(Q)$
may be interpreted as the intervals $\binom{y}{x}_{Q}=[x,y]_{Q}$
in the poset $Q$ and $\preccurlyeq $ as the reversed containment
relation.
\end{defn}

As usual, $\Delta (Q)$ is the order complex of a poset $(Q,\leq
)$. Given a simplicial complex $K$, more generally a polyhedral or
a regular $CW$-complex, the associated face poset is $\Phi
(K)=(\Phi (K),\supseteq )$. Note that $\Phi(K)$ is ordered by the
reversed inclusion, i.e. $F_1\leq F_2$ is equivalent to
$F_1\supseteq F_2$.

\begin{defn}
\label{def_chain_poset} The poset $Chain(Q)=\Phi (\Delta (Q))$ is called the
chain poset associated to $(Q,\leq )$. Its elements are chains $\mathcal{A}%
=\left\{ x_{1}\leq \ldots \leq x_{k}\right\} $ in $Q$ and $\mathcal{%
A\preccurlyeq B}$ if $\mathcal{B}$ is a subchain of $\mathcal{A}$.
\end{defn}

\bigskip

\section{Bredon's theorem}

\bigskip

A fundamental tool in the theory of transformation groups is a
theorem of
Bredon which gives a necessary and sufficient conditions on a $G$-map $%
f:X\longrightarrow Y$ to be a $G$-homotopy equivalence, cf.\
\cite{Bredon} Ch.\ II or \cite{tomDieck} Section II.2. In this paper we need a $\mathbb{Z}%
_{2}$-version of this result. Here and elsewhere throughout the paper we
consistently assume that all spaces are simplicial $\mathbb{Z}_{2}$%
-complexes (polyhedral, $CW$) and that the $\mathbb{Z}_{2}$-maps
are simplicial (cellular). Bredon's theorem holds in higher
generality \cite{JamesSegal} than stated/used in this paper but in
combinatorial applications we can usually restrict our attention
to narrower and more manageable classes of spaces.

\begin{thm}
\label{thm_bredon} \ Suppose that $f:X\longrightarrow Y$ is a
(simplicial) $\mathbb{Z}_{2}$-map of simplicial
$\mathbb{Z}_{2}$-complexes $X$ and $Y$. Let $X^{\mathbb{Z}_{2}}$
and $Y^{\mathbb{Z}_{2}}$ be the associated subspaces of fixed
points and $f^{\mathbb{Z}_{2}}:X^{\mathbb{Z}_{2}}\longrightarrow
Y^{\mathbb{Z}_{2}}$ the map induced by $f$. Then $f$ is
a $\mathbb{Z}_{2}$-homotopy equivalence if and only if both \ $%
f:X\longrightarrow Y$ and $f^{\mathbb{Z}_{2}}:X^{\mathbb{Z}%
_{2}}\longrightarrow Y^{\mathbb{Z}_{2}}$ are homotopy equivalences.
\end{thm}

\begin{cor} If in Theorem~\ref{thm_bredon} the actions of  $\mathbb{Z}_{2}$
on both $X$ and $Y$ are free, i.e.\ if
$X^{\mathbb{Z}_{2}}=Y^{\mathbb{Z}_{2}}=\emptyset$, then a
$\mathbb{Z}_{2}$-map $f : X\longrightarrow Y$  is a
$\mathbb{Z}_{2}$-equivalence if and only it is an ordinary
homotopy equivalence.
\end{cor}

\begin{cor}
\label{cor_bredon} Suppose that $P$ and $Q$ are two involutive posets
(Definition \ref{def_involutive}) and let $f:P\longrightarrow Q$ be a $%
\mathbb{Z}_{2}$-equivariant map of posets. Let $P^{\mathbb{Z}_{2}}$ and $Q^{%
\mathbb{Z}_{2}}$ be the associated subposets of fixed elements. Then $%
f:P\longrightarrow Q$ is a $\mathbb{Z}_{2}$-equivalence if and only if both $%
f:P\longrightarrow Q$ and $f:P^{\mathbb{Z}_{2}}\longrightarrow Q^{\mathbb{Z}%
_{2}}$ are homotopy equivalences of posets.
\end{cor}

Once we reduced the question of $\mathbb{Z}_{2}$-equivalence to
the problem of verifying ordinary homotopy equivalences, we have
on our disposal all the usual combinatorial tools, cf.\
\cite{Bjorner} and \cite{ZZTop}. Our main tool in this paper is
the well known Quillen fiber theorem \cite{Bjorner}
\cite{Quillen} \cite{ZZTop} which says that a monotone map $%
f:P\longrightarrow Q$ of posets is a homotopy equivalence if
$f^{-1}(Q_{\leq q})$ is contractible for each $q\in Q$. Equally
important and useful is the following result widely known as {\em
Order Homotopy Theorem}, \cite{Segal} \cite{Quillen}
\cite{Bjorner}, see also \cite{WeZieZiv} for subsequent
developments and related references.

\begin{prop}
\label{prop_segal} Suppose that $f$ and $g$ are two monotone maps of posets $%
P$ and $Q$ such that $f(x)\leq g(x)$ for each $x\in P.$Then there is a
homotopy equivalence $\Delta (f)\simeq \Delta (g):\Delta (P)\longrightarrow
\Delta (Q)$ between the induced maps of associated order complexes.
Moreover, if $P=Q$ and $g=1_{P}$ $\ $is the identity map, then the
subcomplex ${\rm Im}(f)\subseteq \Delta (P)$ is a deformation retract of $%
\Delta (P)$.
\end{prop}

\begin{exmp}
\label{ex_poznato} The well known fact that the inclusion map $\mathcal{L}%
(P)\longrightarrow P$ is a homotopy equivalence, actually an
inverse to a deformation retraction, is easily deduced from the
second half of Proposition \ref{prop_segal}. Indeed, it is
sufficient to define $f$ as the map $C^{2}:P\longrightarrow P$.
\end{exmp}

\section{$\mathbb{Z}_{2}$-homotopy equivalences of $\mathbb{Z}_{2}$-posets}
\label{sec_all_equivalent}

\bigskip

Suppose that $(P,C)$ is a $WI$-poset (Definition \ref{def_weakly}). Then the
associated Lov\'{a}sz poset $\mathcal{L}(P)$ is involutive with the action
(antitone involution) $\omega :\mathcal{L}(P)\longrightarrow \mathcal{L}(P)$
defined by $\omega (x):=Cx$. The box poset $\mathfrak{B}(P)$ (Definition \ref%
{def_box_poset}) also admits a $\mathbb{Z}_{2}$-action defined by $\omega
(x,y):=(y,x)$. If $(Q,\leq )$ is a $\mathbb{Z}_{2}$-poset with an antitione
involution $\omega :Q\longrightarrow Q$, then both the poset of intervals $%
Int(Q)$ (Definition \ref{def_int_poset}) and the chain poset $Chain(Q)$
(Definition \ref{def_chain_poset}) admit natural $\mathbb{Z}_{2}$-actions.
More precisely, if $\binom{y}{x}\in Int(Q)$ then $\omega \binom{y}{x}:=%
\binom{\omega (x)}{\omega (y)}$ and for $\mathcal{A}=\left\{ x_{1}\leq
\ldots \leq x_{k}\right\} \in Chain(Q),\omega (\mathcal{A})=\mathcal{B}$
where $\mathcal{B}=\left\{ \omega (x_{k})\leq \ldots \leq \omega
(x_{1})\right\} .$

Consequently for each $WI$-poset $(P,C)$ there arise four different $%
\mathbb{Z}_{2}$-posets $\mathcal{L}\mathit{(}P\mathit{)}$, $\mathfrak{B}(P)$%
, $Int(\mathcal{L}(P))$, and $Chain(\mathcal{L}(P)).$ Our
objective is to demonstrate that all these posets are
$\mathbb{Z}_{2}$-homotopy equivalent.

\begin{prop}
\label{prop_omega} Assume that $(Q,\leq )$ is an $I$-poset with an antitone
involution $C:Q\longrightarrow Q$. Then the $\mathbb{Z}_{2}$-map $\Omega
:Int(Q)\longrightarrow \mathfrak{B}(Q)$ defined by $\Omega \binom{y}{x}%
:=(x,Cy)$ is a $\mathbb{Z}_{2}$-isomorphism of $\mathbb{Z}_{2}$-posets.
\end{prop}

\medskip\noindent
{\bf Proof.} Define the inverse map $\Omega ^{\prime
}:\mathfrak{B}(Q)\longrightarrow Int(Q)$ by the formula $\Omega
^{\prime }(a,b):=\binom{Cb}{a}.$ Note that both $\Omega $ and
$\Omega ^{\prime }$ are well defined. It remains to be
shown that one of them, say $\Omega $, is both monotone and $\mathbb{Z}_{2}$%
-equivariant. Indeed, $\Omega $ is monotone since $\binom{y}{x}\preccurlyeq
\binom{y^{\prime }}{x^{\prime }}$ $\iff x\leq x^{\prime }\leq y^{\prime
}\leq y$ implies $x\leq x^{\prime }$ and $Cy\leq Cy^{\prime }$, i.e. $%
(x,Cy)\leq (x^{\prime },Cy^{\prime })$ in $\mathfrak{B}(Q)$. It is $\mathbb{Z%
}_{2}$-equivariant since%
\begin{equation*}
\Omega (\omega \binom{y}{x})=\Omega
\binom{Cx}{Cy}=(Cy,C^{2}x)=(Cy,x)=\omega (x,Cy)=\omega (\Omega
\binom{y}{x}). \quad \blacksquare
\end{equation*}

\begin{prop}
\label{prop_theta} Let $(P,C)$ be a $WI$-poset. Then the
$\mathbb{Z}_{2}$-map\ $\Theta
:\mathfrak{B}(\mathcal{L}(P))\longrightarrow \mathfrak{B}(P)$ of $%
\mathbb{Z}_{2}$-posets, induced by the inclusion map $\mathcal{L}%
(P)\longrightarrow P,$ is a $\mathbb{Z}_{2}$-equivalence.
\end{prop}

\begin{pf} {By Bredon's theorem (Theorem \ref{thm_bredon}), we are supposed
to show that the following two conditions are satisfied,
{\begin{itemize}
\item $\quad \beta :\mathfrak{B}(\mathcal{L}(P))\longrightarrow \mathfrak{B}%
(P)$ is a homotopy equivalence,
\item $\quad \beta ^{\mathbb{Z}_{2}}:\mathfrak{B}(\mathcal{L}(P))^{\mathbb{Z}%
_{2}}\longrightarrow \mathfrak{B}(P)^{\mathbb{Z}_{2}}$ is a
homotopy equivalence.
\end{itemize}}} \par\noindent
Suppose that $(u,v)\in \mathfrak{B}(P).$ Then $(C^{2}u,C^{2}v)\in \mathfrak{B%
}(\mathcal{L}(P))$ and $(C^{2}u,C^{2}v)\geq (u,v)$. Moreover, if $(x,y)\in
\mathfrak{B}(\mathcal{L}(P))$ such that $(x,y)\geq (u,v)$, a consequence of
monotonicity of $C^{2}$ is $(x,y)=(C^{2}x,C^{2}x)\geq (C^{2}u,C^{2}v).$ In
other words $\beta ^{-1}(\mathfrak{B}(P)_{\geq (u,v)})$ has a minimum
element $(C^{2}u,C^{2}v)$, hence it is contractible. By Quillen fiber
theorem $\beta :\mathfrak{B}(\mathcal{L}(P))\longrightarrow \mathfrak{B}(P)$
is a homotopy equivalence.

Let us start with an observation that $\mathfrak{B}(\mathcal{L}(P))^{\mathbb{%
Z}_{2}}=\left\{ (x,x)\mid x\in \mathcal{L}(P)\right\} \cong \mathcal{L}(P)$
and $\mathfrak{B}(P)^{\mathbb{Z}_{2}}=\left\{ (u,u)\mid u\in P\right\} \cong
P.$ It follows from Example~\ref{ex_poznato} that $\beta ^{\mathbb{Z}%
_{2}}$ is also an equivalence of posets.
\end{pf}

\begin{prop}
\label{prop_sigma} Suppose that $(Q,C)$ is an (antitone) involutive poset ($%
I $-poset) in the sense of Definition~\ref{def_involutive}. Let
$\Sigma :Chain(Q)\longrightarrow Int(Q)$ be the map of the
associated chain and interval posets defined by $\Sigma
(\mathcal{A})=\binom{x_{m}}{x_{1}}$ where $\mathcal{A}=\left\{
x_{1}\leq \ldots \leq x_{m}\right\} \in Chain(Q)$. Then $\Sigma $
is a $\mathbb{Z}_{2}$-equivalence.
\end{prop}

\begin{pf}
As before, owing to Bredon's theorem, it is sufficient to show that both $%
\Sigma :Chain(Q)\longrightarrow Int(Q)$ and $\Sigma ^{\mathbb{Z}%
_{2}}:Chain(Q)^{\mathbb{Z}_{2}}\longrightarrow Int(Q)^{\mathbb{Z}_{2}}$ are
homotopy equivalences. Given $\binom{b}{a}\in Int(Q)$ and $\mathcal{A}%
=\left\{ x_{1}\leq \ldots \leq x_{m}\right\} \in Chain(Q)$, we observe that
\begin{equation*}
\binom{b}{a}\leq \Sigma (\mathcal{A})\iff a\leq x_{1}\leq \ldots \leq
x_{m}\leq b.
\end{equation*}
Let $D:=\Sigma ^{-1}(Int(Q)_{\geq \binom{b}{a}})$. Define two monotone maps $%
\lambda ,\mu :D\longrightarrow D$, by the formulas
\begin{equation*}
\lambda (\mathcal{A})=\mathcal{A}^{\prime }:=\left\{ a\leq x_{1}\leq \ldots
\leq x_{m}\leq b\right\} \text{ and }\mu (\mathcal{A}):=\left\{ a\leq
b\right\} .
\end{equation*}%
Let $1_{D}:D\longrightarrow D$ be the identity map. Then $1_{D}(\mathcal{A}%
)\succcurlyeq \lambda (\mathcal{A})\preccurlyeq \mu (\mathcal{A})$ for each $%
\mathcal{A}\in D$. By Proposition \ref{prop_segal} the poset $D$ is
contractible, so by Quillen fiber theorem $\Sigma $ is a homotopy
equivalence. Let us establish now a similar fact for the map $\Sigma ^{%
\mathbb{Z}_{2}}:Chain(Q)^{\mathbb{Z}_{2}}\longrightarrow Int(Q)^{\mathbb{Z}%
_{2}}$. We start with an observation that $\binom{b}{a}\in Int(Q)^{\mathbb{Z}%
_{2}}$ if and only if $\binom{b}{a}=\binom{Cx}{x}$ for some $x\in Q$ such
that $x\leq Cx$. Similarly, $\mathcal{A\in }Chain(Q)^{\mathbb{Z}_{2}}$ if
and only if there exist elements $x_{1}\leq \ldots \leq x_{k}$ in $Q$ such
that $x_{k}\leq Cx_{k}$, in which case $\mathcal{A}:=\left\{ x_{1}\leq
\ldots \leq x_{k}\leq Cx_{k}\leq \ldots \leq Cx_{1}\right\} $. Note that the
inequality $x_{k}\leq Cx_{k}$ is not necessarily strict. Since $\Sigma (%
\mathcal{A})=\binom{Cx_{1}}{x_{1}}$, we observe that $\binom{Ca}{a}\leq
\Sigma (\mathcal{A})$ if and only if $a\leq x_{1}$ and $Cx_{1}\leq Ca$.
Define $\lambda _{1}(\mathcal{A})=\mathcal{A}^{\prime }:=\left\{ a\leq
x_{1}\leq \ldots \leq x_{k}\leq Cx_{k}\leq \ldots \leq Cx_{1}\leq Ca\right\}
$ and $\mu _{1}(\mathcal{A}):=\left\{ a\leq Ca\right\} $ as monotone maps on
the poset $D_{1}:=(\Sigma ^{\mathbb{Z}_{2}})^{-1}(Int(Q)_{\geq \binom{Ca}{a}%
}^{\mathbb{Z}_{2}})$. Since $1_{D_{1}}(\mathcal{A})\succcurlyeq \lambda _{1}(%
\mathcal{A})\preccurlyeq \mu _{1}(\mathcal{A})$ for each $\mathcal{A}\in
D_{1}$ we deduce from Proposition \ref{prop_segal} that $D_{1}$ is
contractible. Hence, by Quillen fiber theorem, $\Sigma ^{\mathbb{Z}_{2}}$ is
a homotopy equivalence which completes the proof of the proposition.
\end{pf}

\begin{cor}
\label{cor_sigma} Let $(P,C)$ be a $WI$-poset and $\mathcal{L}(P)$ the
associated Lov\'{a}sz subposet. Then $\Sigma :Chain(\mathcal{L}%
(P))\longrightarrow Int(\mathcal{L}(P))$ is a $\mathbb{Z}_{2}$-homotopy
equivalence.
\end{cor}

\begin{prop}
\label{prop_gama} Suppose that $(Q,C)$ is an involutive poset ($I$-poset).
Then there is a $\mathbb{Z}_{2}$-homotopy equivalence of $\mathbb{Z}_{2}$%
-complexes $\Delta (Q)$ and $\Delta (Chain(Q))$.
\end{prop}

\begin{pf}
Note that $\Delta (Chain(Q))$ is just the first baricentric subdivision of
the simplicial complex $\Delta (Q)$. Hence, there is a well known canonical
homeomorphism $\Gamma :|\Delta (Q)|\longrightarrow |\Delta (Chain(Q))|$ \ of
the associated geometric realizations of these complexes. What remains to be
done is to show that $\Gamma $ is $\mathbb{Z}_{2}$-equivariant. Recall that
the $\mathbb{Z}_{2}$-actions on $Q$ and $Chain(Q)$ are given by $\omega
(q)=Cq$ and $\omega \left\{ x_{1}\leq \ldots \leq x_{k}\right\} =\left\{
Cx_{k}\leq \ldots \leq Cx_{1}\right\} $ respectively. The homeomorphism $%
\Gamma $ is explicitly defined as follows. Let
$t=t_{1}x_{1}+\ldots +t_{k}x_{k}\in |\Delta (Q)|$, where
$x_{1}<\ldots <x_{k}$, $t_{1}+\ldots +t_{k}=1$, and $t_{j}\geq 0$.
Put the sequence $(t_{j})_{j=1}^{k}$in the descending order which
means that for some permutation $\pi :[n]\longrightarrow \lbrack
n]$ we have inequalities $t_{\pi _{1}}\geq t_{\pi _{2}}\geq \ldots
\geq t_{\pi _{k}}$. Then $X_{\pi _{1}}\succcurlyeq X_{\pi
_{2}}\succcurlyeq \ldots \succcurlyeq X_{\pi _{k}}$, where $X_{\pi
_{j}}$ :=$\left\{ x_{\pi _{1}},x_{\pi _{2}},\ldots ,x_{\pi
_{j}}\right\} $, is a chain in the poset $Chain(Q)$ and $\Gamma
(t)=s_{1}X_{\pi _{1}}+s_{2}X_{\pi _{2}}+\ldots +s_{k}X_{\pi
_{k}}\in |\Delta (Chain(Q))|$,
where the relation between sequences $(s_{j})_{j=1}^{k}$ and $%
(t_{j})_{j=1}^{k}$ is determined by the following equality
\begin{equation*}
s_{1}x_{\pi _{1}}+s_{2}\frac{x_{\pi _{1}}+x_{\pi _{2}}}{2}+\ldots +s_{k}%
\frac{x_{\pi _{1}}+\ldots +x_{\pi _{k}}}{k}=t_{1}x_{1}+\ldots +t_{k}x_{k}%
\text{.}
\end{equation*}%
Note that \ $\omega (t)=t_{1}C(x_{1})+\ldots +t_{k}C(x_{k})$, and $\Gamma
\omega (t)=s_{1}(\omega X_{\pi _{1}})+s_{2}(\omega X_{\pi _{2}})+\ldots
+s_{k}(\omega X_{\pi _{k}})=\omega (\Gamma (t))$ which implies that $\Gamma $
is indeed $\mathbb{Z}_{2}$-equivariant.
\end{pf}

\bigskip

All results in this section together imply that there exists essentially a
unique $\mathbb{Z}_{2}$-homotopy type associated to a given $WI$-poset $%
(P,C) $.

\begin{cor}
For a $WI$-poset $(P,C)$, the order complexes of $\mathbb{Z}_{2}$%
-posets
\begin{equation*}
\mathcal{L}(P)\mathit{\quad }\text{ \ }\mathfrak{B}(P)\quad \text{ \ \ }%
\mathfrak{B}(\mathcal{L}(P))\quad Int(\mathcal{L}(P))\text{ \quad\ }Chain(%
\mathcal{L}(P))
\end{equation*}%
are all $\mathbb{Z}_{2}$-homotopy equivalent.
\end{cor}

\bigskip

\section{Relatives of the box poset }\label{sec_avatars}

\bigskip

As a variation on a theme, motivated by applications in Section \ref%
{sec_applications}, we introduce two more relatives of the box poset $%
\mathfrak{B}(P)$.

\begin{defn}
\label{def_extended_box} Assume that $(P,C)$ is a $WI$-poset. Define $%
\widehat{P}:=P\cup \{ \widehat{0}\} $ as a new poset
obtained by adding to $P$ a possibly new minimum element $\widehat{0}$. The \emph{%
extended box poset} $\mathfrak{B}_{ex}(P)$, associated to the $WI$-poset $P,$
is a subposet of $\widehat{P}\times \widehat{P}$ defined by $\mathfrak{B}%
_{ex}(P):=\mathfrak{B}(P)\cup \{ (p,\widehat{0})\mid p\in P\} \cup
\{ (\widehat{0},q)\mid q\in P\} $.
\end{defn}

\begin{thm}
\label{thm_extended_box} Suppose that $(P,C)$ is a $WI$-poset and $\mathfrak{%
B}(P)$, $\mathfrak{B}_{ex}(P)$ the box poset, respectively the extended box
poset associated to $P$. Then the inclusion map $e:\mathfrak{B}%
(P)\longrightarrow \mathfrak{B}_{ex}(P)$ is a $\mathbb{Z}_{2}$-equivalence
of posets.
\end{thm}

\begin{pf}
As before, we ought to show that both $e:\mathfrak{B}(P)\longrightarrow
\mathfrak{B}_{ex}(P)$ and $e^{\mathbb{Z}_{2}}:\mathfrak{B}(P)^{\mathbb{Z}%
_{2}}\longrightarrow \mathfrak{B}_{ex}(P)^{\mathbb{Z}_{2}}$ are homotopy
equivalences of posets. Let us show that for each $(p,q)\in \mathfrak{B}%
_{ex}(P)$ the poset $D_{p,q}:=e^{-1}(\mathfrak{B}_{ex}(P)_{\geq
(p,q)})$ is contractible. This is obvious if $p\neq
\widehat{0}\neq q$ since in that case
$D_{p,q}=\mathfrak{B}(P)_{\geq (p,q)}$. Let us establish the
contractibility of $D_{p,\hat{0}}$, the case of $D_{\hat{0},q}$ is
treated similarly. By definition $(x,y)\in D_{p,\hat{0}}$ if and
only if $x\geq p$, $x\leq Cy$, and $y\leq Cx$. As a consequence\
we have the inequalities $y\leq Cx\leq Cp$. This means that
$(p,y)\in \mathfrak{B}(P)$
and, since $(p,y)\geq (p,\widehat{0})$, we conclude that $(p,y)\in D_{p,%
\hat{0}}$. The identity map $1_{D_{p,\hat{0}}}$ and the map $\mu
:D_{p,\hat{0}}\longrightarrow D_{p,\hat{0}}$ defined by $\mu
(x,y)=(p,y)$ satisfy the condition $\mu (x,y)=(p,y)\leq (x,y)=1_{D_{p,%
\hat{0}}}(x,y)$ hence, by Proposition \ref{prop_segal}, $E_{p,\hat{0}%
}={\rm Im}(\mu )=\left\{ (p,y)\mid y\leq Cp\And p\leq Cy\right\} $
is a deformation retract of $D_{p,\hat{0}}$. On the other hand,
since $(p,Cp)$
is the maximum element of $E_{p,\hat{0}}$, we conclude that $E_{p,%
\hat{0}}$ is contractible, so the same holds for $D_{p,\hat{0}}$.

The case of the map $e^{\mathbb{Z}_{2}}:\mathfrak{B}(P)^{\mathbb{Z}%
_{2}}\longrightarrow \mathfrak{B}_{ex}(P)^{\mathbb{Z}_{2}}$ is simpler since
$\mathfrak{B}(P)^{\mathbb{Z}_{2}}=\mathfrak{B}_{ex}(P)^{\mathbb{Z}%
_{2}}=\left\{ (p,p)\mid p\in P\right\} $ and $e^{\mathbb{Z}_{2}}$ is an
identity map.
\end{pf}

\bigskip

\begin{defn}
\label{def_enriched_box} \ Suppose that $(P,C)$ is a $WI$-poset
and that $P$ is a subposet of an auxiliary poset $S$. Define the
box poset of $P$ enriched over $S$ as the
$\mathbb{Z}_{2}$-subposet of $\widehat{S}\times \widehat{S}$
described by the equality
\begin{equation*}
\mathfrak{B}_{S}(P):=\mathfrak{B}(P)\cup \left\{ (p,\widehat{0})\mid p\in
S\right\} \cup \left\{ (\widehat{0},q)\mid q\in S\right\} .
\end{equation*}
\end{defn}

\begin{thm}
\label{thm_enriched_box} Suppose that $(P,C)$ is a $WI$-poset, $S$ a
superposet of $P,$ and $\mathfrak{B}_{S}(P)$ the box poset of $P$ enriched
over $S$. If $S$ is contractible then the geometric realization $|\Delta (%
\mathfrak{B}_{S}(P))|$ of this poset is a $\mathbb{Z}_{2}$-space which is $%
\mathbb{Z}_{2}$-homotopy equivalent to the suspension $Susp(|\Delta (%
\mathfrak{B}(P))|)$ of the geometric realization of the box poset $\mathfrak{%
B}(P)$.
\end{thm}

\begin{pf}
Let $\overleftrightarrow{\mathfrak{B}(P)}=\mathfrak{B}(P)\cup \left\{
a_{1},a_{2}\right\} $ be the poset obtained from the box poset $\mathfrak{B}%
(P)$ by adding two new incomparable minimal elements $a_{1}$ and $a_{2}$.
Extend the involution $\omega $ from $\mathfrak{B}(P)$ to $%
\overleftrightarrow{\mathfrak{B}(P)}$ by the requirement that $\omega
(a_{1})=a_{2}$ $\ $and $\omega (a_{2})=a_{1}$. Note that each chain $%
\mathcal{A}$ in $\mathfrak{B}(P)$ can be extended to chains $\mathcal{A}_{1}=%
\mathcal{A\cup }\left\{ a_{1}\right\} $ and $\mathcal{A}_{2}=\mathcal{A\cup }%
\left\{ a_{2}\right\} $. Since $a_{1}$ and $a_{2}$ are incomparable, we
observe that the geometric realization of the order complex $\Delta (%
\overleftrightarrow{\mathfrak{B}(P)})$ is, as a $\mathbb{Z}_{2}$-space,
homeomorphic to the suspension $Susp(|\Delta (\mathfrak{B}(P))|)$. Define a
monotone, $\mathbb{Z}_{2}$-map $\Psi :\mathfrak{B}_{S}(P)\longrightarrow
\overleftrightarrow{\mathfrak{B}(P)}$ of posets as follows. If $p,q\in P$
then $\Psi (p,q)=(p,q).$ Otherwise $\Psi (p,\widehat{0})=a_{1}$ and $\Psi (%
\widehat{0},p)=a_{2}$ for each $p\in S$. The map $\Psi $ is obviously $%
\mathbb{Z}_{2}$-equivariant. Let us show that it is a $\mathbb{Z}_{2}$%
-equivalence. In light of Theorem \ref{thm_bredon} we ought to show that
both $\Psi $ and $\Psi ^{\mathbb{Z}_{2}}$ are ordinary homotopy
equivalences. The map $\Psi ^{\mathbb{Z}_{2}}$ turns out to be essentially
an identity map so we focus our attention on $\Psi $. Let $D_{x,y}:=\Psi
^{-1}(\overleftrightarrow{\mathfrak{B}(P)}_{\geq (x,y)})$. If $(x,y)\in
\mathfrak{B}(P)$ then D$_{x,y}=\mathfrak{B}_{S}(P)_{\geq (x,y)}$, hence it
is contractible. If $(x,y)=(p,\widehat{0})$ for some $p\in S$, then $D_{p,%
\hat{0}}$ can be deformed to its subposet $S_{1}:=\{ (p,\widehat{0}%
)\mid p\in S\} .$ Indeed, such a deformation is provided by the map $%
\mu :D_{p,\hat{0}}\longrightarrow S_{1}$, where $\mu (x,y):=(x,\widehat{0%
})$. Since $S_{1}\cong S$, and by assumption $S$ is contractible, we
conclude that $D_{p,\hat{0}}$ is contractible. By a similar argument $D_{%
\hat{0},q}$ is also contractible and finally, by Quillen fiber theorem, $%
\Psi $ is a homotopy equivalence of posets.
\end{pf}

\bigskip

\section{Applications to graph complexes\label{sec_applications}}

\bigskip

Suppose that $G=(V_{G},E_{G})$ is a finite graph. The poset
$(P_{G},\subseteq )$, where by definition $P_{G}:=\left\{ A\subset
V_{G}\mid CN(A)\neq \varnothing \right\} $, is weakly involutive
(a $WI$-poset) where the weak involution $C:P_{G}\longrightarrow
P_{G}$ is defined by $C(A$)$:=CN(A).$
This is precisely the example which served as a motivation for introducing $%
WI$-posets and the development of the associated $\mathbb{Z}_{2}$-posets ($%
\mathbb{Z}_{2}$-complexes). By specialization, each of the $\mathbb{Z}_{2}$%
-posets from Sections \ref{sec_all_equivalent} and
\ref{sec_avatars} yields the corresponding graph
$\mathbb{Z}_{2}$-complex. Here is a partial list of these
complexes
\begin{equation*}
\mathcal{L}(G):=\mathcal{L}(P_{G}) \quad  \mathfrak{B}(G):=
\mathfrak{B}(P_{G}) \quad
\mathfrak{B}_{ex}(G):=\mathfrak{B}(P_{G})\quad
\mathfrak{B}_{S}(G):=\mathfrak{B}_{S}(P_{G}).
\end{equation*}
One of our objectives in this section is to compare these
complexes with the existing graph complexes listed in
\cite{MatZieg}. More importantly, we demonstrate that in virtually
all cases analyzed in \cite{MatZieg} (Theorems
1 and 3), the equality $Ind_{\mathbb{Z}_{2}}(K_{1})=Ind_{\mathbb{Z}%
_{2}}(K_{2})$ of $\mathbb{Z}_2$-indices of graph complexes is a
consequence of
the stronger statement that $K_{1}$ and $K_{2}$ are $\mathbb{Z}_{2}$%
-homotopy equivalent. Similarly, all the inequalities $Ind_{\mathbb{Z}%
_{2}}(K_{1})\leq Ind_{\mathbb{Z}_{2}}(K_{1})+1$ are found to be consequences
of the\ $\mathbb{Z}_{2}$-equivalence $K_{1}\cong Susp(K_{2})$.

\medskip\noindent
{\bf Caveat} We interchangeably use the words posets and complexes
for the same objects. This should not cause any ambiguity since
one already talks about the homology and the homotopy of a poset
$Q$, having in mind the homology and homotopy of the associated
order complex $\Delta(Q)$.

\medskip\noindent
The complex $\mathcal{L}(G)$ is of course the Lov\'{a}sz original $\mathbb{Z}%
_{2}$-poset ($\mathbb{Z}_{2}$-complex), denoted by $L(G)$ in \cite{MatZieg}%
. The poset \ $\mathfrak{B}(G)$ is easily identified as the box complex $%
B_{chain}(G)$, while $\mathfrak{B}_{ex}(G)$ is clearly the box complex $B(G)$
from \cite{MatZieg}. The complex $B_{0}(G)$ is recognized as our complex $%
\mathfrak{B}_{S}(G)$ where $S=$ $\mathcal{P}^{\prime }(V_{G}):=\mathcal{P}%
(V_{G})\smallsetminus \left\{ \varnothing \right\} $ is the poset
of all non-empty subsets of $V_{G}$.

All $\mathbb{Z}_{2}$-equivalences between these complexes (and
their suspensions) are immediate consequences of results from Sections \ref%
{sec_all_equivalent} and \ref{sec_avatars}. The complexes from \cite{MatZieg}
that do not automatically fit into this scheme are complexes $B_{edge}(G)$, $%
B_{Sark}^{KG}(\mathcal{F})$, $B_{chain}^{KG}(\mathcal{F})$, listed
as complexes no. 4, 5, and 6 in Section 5 of \cite{MatZieg}. Note
that the complexes {\[ B_{Sark}^{KG}(\mathcal{F}):=\Delta \{
B^{\prime }\uplus B^{\prime \prime }\mid B^{\prime },\, B^{\prime
\prime }\subseteq \lbrack n],\, B^{\prime
}\cap B^{\prime \prime }=\varnothing ,\,(\exists X\in \mathcal{F}) X\subseteq A%
\text{ or }X\subseteq B \} \]
\[
 B_{chain}^{KG}(\mathcal{F}):=\Delta \{ B^{\prime
}\uplus B^{\prime \prime }\mid B^{\prime },\, B^{\prime \prime
}\subseteq \lbrack n],\, B^{\prime }\cap B^{\prime \prime
}=\varnothing , \, (\exists X,Y\in \mathcal{F})X\subseteq A\text{
and }Y\subseteq B\}
\]}
\noindent are defined in terms of the chosen Kneser representative
$\mathcal{F}$ of the graph $G=KG(\mathcal{F})$. This explains why
they cannot be immediately expressed in terms of the associated
$WI$-posets $(P_{G},C)$. Nevertheless, the approach based on
Bredon's theorem (Theorem \ref{thm_bredon}) is equally efficient
and elegant.

\begin{prop}
\label{prop_chain} Let $G=KG(\mathcal{F})$ be the Kneser graph associated to
a finite family of sets $\mathcal{F}$. Let $\mathfrak{B}%
(G)=B_{chain}(G) $ and $B_{chain}^{KG}(\mathcal{F})$ be the
associated box complexes {\rm (}no.\ {\rm 3} and no.\ {\rm 6}
from the list in Section~5 of {\rm \cite{MatZieg}}%
{\rm )}. Then the map%
\begin{equation*}
\Phi :B_{chain}(G)\longrightarrow B_{chain}^{KG}(\mathcal{F})
\end{equation*}%
defined by $\Phi (a\uplus b):=(A\uplus B)$ where $A:=\cup a$ and $B:=\cup b$%
, is a $\mathbb{Z}_{2}$-homotopy equivalence of $\mathbb{Z}_{2}$-posets {\rm (}$%
\mathbb{Z}_{2}$-spaces{\rm )}.
\end{prop}

\begin{pf}
The action of $\mathbb{Z}_{2}$ on both $B_{chain}^{KG}(G)$ and $%
B_{chain}^{KG}(\mathcal{F})$ is free hence, in light of Theorem \ref%
{thm_bredon}, it is sufficient to show that $\Phi $ is a homotopy
equivalence. Again, the Quillen fiber theorem proves to be a very
convenient tool. Given $A\uplus B\in B_{chain}^{KG}(\mathcal{F})$,
let
\begin{equation*}
D_{A,B}:=\Phi ^{-1}(B_{chain}^{KG}(\mathcal{F})_{\geq A\uplus B}):=\left\{
a\uplus b\in B_{chain}^{KG}(G)\mid \cup a\subseteq A\text{ and }\cup
b\subseteq B\right\} .
\end{equation*}%
Note that both $a^{\prime }:=\left\{ X\in \mathcal{F}\mid
X\subseteq A\right\} $ and \ $b^{\prime }:=\left\{ Y\in
\mathcal{F}\mid Y\subseteq B\right\} $ are non-empty. Moreover, $\
a^{\prime }\uplus b^{\prime }$ is the maximum element in
$D_{A,B}$, hence $D_{A,B}$ is contractible. It immediately follows
that $\Phi $ is a homotopy equivalence and, a posteriori by
Bredon's theorem, $\Phi $ is a $\mathbb{Z}_{2}$-homotopy
equivalence.
\end{pf}

\begin{prop}
\label{prop_sark} Assume that $G=KG(\mathcal{F})$ is the Kneser graph
associated to $\mathcal{F}$ and let $B_{0}(G)\cong \mathfrak{B}_{\mathcal{P}%
^{\prime }(V)}(G)$ and $B_{Sark}^{KG}(\mathcal{F})$ be the box
complexes (posets) no.\ {\rm 2} and no.\ {\rm 6} from the list in
Section 5 of \cite{MatZieg}. Then the map
\begin{equation*}
\Psi :B_0(G)\longrightarrow B_{Sark}^{KG}(\mathcal{F})
\end{equation*}%
defined by $\Psi (a\uplus b):=(A\uplus B)$ is a $\mathbb{Z}_{2}$-homotopy
equivalence.
\end{prop}

\begin{pf}
The proof is similar to the proof of Proposition \ref{prop_chain}. If $%
A\uplus B\in B_{Sark}^{KG}(\mathcal{F})$, then by definition at least one of
the sets $A$ and $B$ contains an element $X\in \mathcal{F}$ as a subset. If
both $A$ and $B$ satisfy this condition then, as in the proof of Proposition %
\ref{prop_chain}, the set $D_{A,B}:=\Psi ^{-1}(B_{Sark}^{KG}(\mathcal{F}%
)_{\geq A\uplus B})$ has a maximum element and must be contractible. Suppose
that $A\supseteq X\in \mathcal{F}$ but $B$ does not contain elements from $%
\mathcal{F}$ as subsets. Let $c:=\left\{ Y\in \mathcal{F}\mid Y\subseteq
A\right\} $. Define $\mu :D_{A,B}\longrightarrow D_{A,B}$ as the monotone
map such that $\mu (x\uplus y):=(c\uplus y)$. Since always $x\uplus y\geq
c\uplus y$, we conclude that ${\rm Im}(\mu )$ is a deformation retract of $%
D_{A,B}$. On the other hand \ ${\rm Im}(\mu )$ has the maximum element \ $%
c\uplus \varnothing $, hence it is contractible. This again allows us to use
Quillen fiber theorem to conclude that $\Psi $ is a homotopy equivalence.
Bredon's theorem as before implies that $\Psi $ is actually a $\mathbb{Z}%
_{2} $-homotopy equivalence.
\end{pf}

For completeness we formulate one more result involving the
complex $B_{edge}(G)$, listed as no. 4 in the list in Section 5 of
\cite{MatZieg}. Recall that  \[ B_{edge}(G):=\{F\subset A'\times
A'' \vert \emptyset\neq A', A''\subset V, \, A'\cap A''=\emptyset,
\, G[A',A''] \text{ is complete}\}. \]
\begin{prop}
Let $G=KG(\mathcal{F})$ be the Kneser graph associated to a finite
family of sets $\mathcal{F}$. Let $\mathfrak{B}(G)=B_{chain}(G)$
and $B_{edge}(G)$ be the associated
box complexes {\rm (}no. 3 and no. 4 from the list in Section~5 of {\rm \cite{MatZieg}}%
{\rm )}. Then the map%
\begin{equation*}
\Lambda :B_{edge}(G)\longrightarrow B_{chain}(G)
\end{equation*}%
defined by $\Lambda (F):=(A'\uplus A'')$ for $F\subset A'\times
A''$, is a $\mathbb{Z}_{2}$-homotopy equivalence of
$\mathbb{Z}_{2}$-posets.
\end{prop}

\begin{pf}
The proof is similar in spirit to the proofs in this and earlier
sections so the details are omitted.
\end{pf}

We have convinced ourselves that all $\mathbb{Z}_2$-complexes
\begin{equation}
\label{eq_lambda} \mathcal{L}(P_G)\quad \mathfrak{B}(P_G)\quad
\mathfrak{B}_{ex}(P_G)\quad B_{chain}(G)\quad B_{edge}(G)
\end{equation}
have the same $\mathbb{Z}_2$-homotopy type.

\begin{defn}\label{def_lambda}
Given a graph $G=(V_G,E_G)$, let $\Lambda(G)$ the common
$\mathbb{Z}_2$-homotopy type of each of the complexes listed in
{\rm (\ref{eq_lambda})}. We occasionally, by a slight abuse of
language, refer to $\Lambda(G)$ as to the graph complex associated
to $G$.
\end{defn}

\bigskip

\section{Which $\mathbb{Z}_2$-complexes are graph complexes?}
\label{sec_which}

\bigskip

Suppose that $\hat{L} = L\cup\{\hat{0},\hat{1}\}$ is a finite
lattice with $L$ as its {\em proper part}, \cite{Bjorner}
Section~3. The proper part $L$ of $\hat{L}$ is a semilattice in
the sense that each subset $A\subset L$, bounded from above, has a
least upper bound, similarly each $B\subset L$ bounded from below
has a greatest lower bound. Conversely, each semilattice $L$ is
the proper part of the lattice $\hat{L}:=
L\cup\{\hat{0},\hat{1}\}$ where $\hat{0}$ and $\hat{1}$ are added
minimum (maximum) elements. Assume that $(L,\leq)$ is a
semilattice which is also an $I$-poset
(Definition~\ref{def_involutive}) with a monotone involution $C:
L\rightarrow L$. Note that $\hat{L} = L\cup\{\hat{0},\hat{1}\}$ is
also an $I$-poset where $C : \hat{L}\rightarrow \hat{L}$ is an
extension the old involution, $C(\hat{0})=0, C(\hat{1})=\hat{1}$.

If $L$ is free in the sense that $C$ satisfies an additional
condition, $x\leq C(x) \Rightarrow x=\hat{0} \text{\rm { or }}
x=\hat{1}$, or equivalently if the involution $C: L\rightarrow L$
is fixed-point-free, then we call $L$ a free $I$-semilattice. The
ortholattices used by Walker, \cite{Walker} Section~4, are very
similar to our free $I$-semilattices, the main difference being
that the involution $C:L\rightarrow L$ in an ortholattice is
antitone, rather than monotone. The condition $x\leq y \Rightarrow
C(y)\leq C(x)$ implies that the ``orthogonality relation'',
$x\perp y \Leftrightarrow x\leq C(y)$, is symmetric which leads to
an ``orthogonality graph'' $G^{\perp}=(V^{\perp}, E^{\perp})$
associated to $L$ defined by $V^{\perp}=L$ and $(x,y)\in
E^{\perp}\Leftrightarrow x\perp y$. The associated (neighborhood)
graph complex (lattice) turns out to be closely related to the
original ortholattice $L$ and among the consequences is the result
that each ortholattice arises as the graph complex (lattice) of
some graph.

Each free $I$-semilattice also can be associated a natural graph
$G_L = (V_L, E_L)$ and our main objective in this section is to
analyze its graph complex $\Lambda(G_L)$.

\begin{defn}
\label{def_comp} Suppose that $(L,\leq)$ is a free $I$-semilattice
i.e.\ a semilattice which is a monotone $I$-poset with a
fixed-point-free involution $C:L\rightarrow L$. Define the
associated ``compatibility graph'' $G_L = (V_L, E_L)$ as the graph
on the ground set $V_L:= L$ such that $(x,y)\in E_L
\Leftrightarrow y\leq C(x) \text{\rm { or }} x\leq C(y)$.
\end{defn}

The {\em ``fat''} semilattices or $F$-semilattices for short, are
particularly well behaved and admit a short and transparent
description of its ``compatibility graph'' $G_L$.
\begin{defn}
\label{def_fat} A semilattice $(L,\leq)$ is a $F$-semilattice if
its intervals $[x,y]_L$ are ``fat'' in the sense that for each
strict chain $x<z<y$ in $L$ there is an element $z'\in [x,y]_L$,
incomparable to $z$.
\end{defn}

Suppose from here on that $(L,\leq)$ is a free $I$-semilattice
with ``fat'' intervals. Let $N(G_L)$ be the neighborhood complex
of $G_L$ and $\mathcal{L}(G_L)$ the associated Lov\'{a}sz complex.
By definition
\[
N(\{x\})= L_{\geq C(x)}\cup L_{\leq C(x)} = C(L_{\geq x}\cup
L_{\leq x}) = C(Comp(\{x\}))
\]
where $Comp(B)$ is the set of all elements in $L$ which are
$\leq$-comparable with all elements $y\in B$. Let us observe that
for each $A\subset L$, if $N(A)\neq\emptyset$ then there exists a
chain $a_1\leq b_1\leq a_2\leq b_2\leq \ldots\leq a_k\leq b_k$ in
$\hat{L}$ such that
\begin{equation}
\label{eqn_intervali} N(A) = [a_1,b_1]_L\cup
[a_2,b_2]_L\cup\ldots\cup [a_k,b_k]_L.
\end{equation}
Note that we allow elements in this chain to be $\hat{0}$ or
$\hat{1}$, however the intervals are always taken in $L$ so for
example $[\hat{0},x]_L = L_{\leq x}$ and $[y,\hat{1}]_L = L_{\geq
y}$. The observation follows by an easy induction on the size of
$A\subset L$. Indeed, if $A' = A\cup\{x'\}$, then $N(A')=N(A)\cap
N(\{x'\}) = N(A)\cap (L_{\geq C(x')}\cup L_{\leq C(x')})$ and, if
$N(A)$ admits a decomposition (\ref{eqn_intervali}), it is easily
checked that $N(A')$ also admits such a decomposition. Note that
here we did not use the fact that $L$ has ``fat'' intervals. This
hypothesis is essentially used in the proof of the following
lemma.

\begin{lemma}
\label{lema_fat} If $N(A)$ admits the decomposition {\rm
(\ref{eqn_intervali})} then
\begin{equation}
\label{eqn_dual} N(N(A)) = [\hat{0},C(a_1)]_L\cup
[C(b_1),C(a_2)]_L\cup \ldots\cup [C(b_{k-1}),C(a_k)]_L\cup
[C(b_k),\hat{1}].
\end{equation}
\end{lemma}

\begin{pf}
Let $Comp(N(A))$ be the set of all elements in $L$ comparable to
all elements in $N(A)$. Since the intervals in $L$ are ``fat'', we
observe that
\[
Comp(N(A)) = [\hat{0},a_1]_L\cup [b_1,a_2]_L\cup \ldots\cup
[b_{k-1},a_k]_L\cup [b_k,\hat{1}]
\]
and the Lemma is deduced from the fact that $N(N(A)) =
C(Comp(N(A))$.
\end{pf}

\begin{thm}
\label{thm_teorema} Suppose that $(L,\leq)$ is a free
$I$-semilattice with ``fat'' intervals, $F$-se\-mi\-latti\-ces in
the sense of Definition~\ref{def_fat}, and let $G_L=(V_L,E_L)$ be
the associated ``compatibility graph''. Then the graph complex
$\Lambda(G_L)$ (Definition~\ref{def_lambda}) of $G_L$ is
$\mathbb{Z}_2$-homotopy equivalent to the order complex
$\Delta(L)$.
\end{thm}

\begin{pf}
Let $\hat{\mathcal{L}}(G_L)= \{ \diamondsuit \subset L \mid
N(N(\diamondsuit)) = \diamondsuit\}$ be the Lov\'{a}sz lattice and
${\mathcal{L}}(G_L)=\hat{\mathcal{L}}(G_L)\setminus\{\hat{0},\hat{1}\}$
its proper part. We already know that $\diamondsuit\in
{\mathcal{L}(G_L)}$ if and only if
$L\neq\diamondsuit\neq\emptyset$ and $\diamondsuit = N(A)$ admits
a decomposition into a union of $L$-intervals, described in
(\ref{eqn_intervali}). The poset $\mathcal{L}(G_L)$ is ordered by
the reversed inclusion, i.e.\ $\diamondsuit_1\leq \diamondsuit_2
\Leftrightarrow \diamondsuit_1 \supseteq \diamondsuit_2$. Our
objective is to compare the semilattice ${\mathcal{L}}(G_L)$ and
the original semilattice $L$. Let $Chain(L)=\Phi(\Delta(L))$ be
the chain poset associated to $(L,\leq)$,
Definition~\ref{def_chain_poset}. By Proposition~\ref{prop_gama},
posets $L$ and $Chain(L)$ are $\mathbb{Z}_2$-homotopy equivalent.
Define the map $\Omega : \mathcal{L}(G_L)\rightarrow Chain(L)$ by
the formula
\[
\Omega([a_1,b_1]_L\cup [a_2,b_2]_L\cup\ldots\cup [a_k,b_k]_L) =
(a_1\leq b_1\leq \ldots \leq a_k\leq b_k)\in Chain(L).
\]
Here as before, some of the elements are allowed to be $\hat{0}$
or $\hat{1}$ but in the chain itself they are neglected.
Alternatively, one can agree that each chain in $L$ is enriched by
elements $\hat{0}$ and $\hat{1}$. The map $\Omega$ is
$\mathbb{Z}_2$-equivariant in the sense that for each
$\diamondsuit\in \mathcal{L}(G_L), \, C\Omega(\diamondsuit) =
\Omega C(\diamondsuit)$. Unfortunately the map $\Omega$ is not
monotone (antitone). This is not a surprise since
$\mathcal{L}(G_L)$ is an antitone while $Chain(L)$ is a monotone
$I$-poset, hence there does not exists a
$\mathbb{Z}_2$-equivariant monotone (antitone) map of these
posets. In order get around this difficulty we pass to the
$I$-poset $Chain(\mathcal{L}(G_L))$ which is a monotone $I$-poset
and which, according to Proposition~\ref{prop_gama}, retains the
$\mathbb{Z}_2$-homotopy type of the poset $\mathcal{L}(G_L)$. The
map $\Omega$ can be extended to a $\mathbb{Z}_2$-equivariant,
monotone map $\Omega^{\sharp} : Chain(\mathcal{L}(G_L))\rightarrow
Chain(L)$ of posets as follows. Given a chain
$\diamondsuit_1\leq\diamondsuit_2\leq\ldots\leq\diamondsuit_k$ in
$\mathcal{L}(G_L)$, the associated elements
$\Omega(\diamondsuit_1), \Omega(\diamondsuit_2),\ldots,
\Omega(\diamondsuit_k)$ are not necessarily elements of a chain in
Chain(L). The obstacle is that they may not be comparable.
However, their union is a well defined chain in $L$ so by
definition
\[
\Omega^{\sharp}((\diamondsuit_j)_{j=1}^k) := \bigcup_{j=1}^k
\Omega(\diamondsuit_j).
\]
We claim that $\Omega^{\sharp}$ is a $\mathbb{Z}_2$-homotopy
equivalence. By Theorem~\ref{thm_bredon}, it is sufficient to show
that $\Omega^{\sharp}$ is an ordinary homotopy equivalence. As
before, the Quillen fiber theorem is a convenient tool. Given a
chain $\Gamma=(c_1\leq c_2\leq\ldots \leq c_k)\in Chain(L)$, our
objective is to show that $D_\Gamma:=
(\Omega^{\sharp})^{-1}(Chain(L)_{\geq \Gamma})$ is a contractible
subposet of $Chain(\mathcal{L}(G_L))$. Note that $D_\Gamma$ is
itself a chain poset, $D_\Gamma = Chain(E_\Gamma)$. Indeed,
$E_\Gamma$ is characterized by the condition $\diamondsuit\in
E_\Gamma$ if and only if $\Omega(\diamondsuit)$ is a subchain of
$\Gamma$. Since $D_\Gamma$ and $E_\Gamma$ have the same homotopy
type, it is sufficient to show that $E_\Gamma$ is contractible.
One way to establish this fact is to observe that $E_\Gamma$ is a
semilattice and that $\diamondsuit_0:=[c_1,c_k]$ is an element in
$E_\Gamma$ with an empty set of complements, cf.\ \cite{Bjorner}
Theorem~10.15. One can also note that $E_\Gamma$ is isomorphic to
the poset $F_\Gamma$ where $I\in F_\Gamma$ if $I = \diamondsuit
\cap \Gamma$ for some $\diamondsuit \in E_\Gamma$. In other words
elements of $F_\Gamma$ are unions of intervals in $\Gamma$. So
there are alternative proofs that $D_\Gamma\cong F_\Gamma$ is
contractible, for example one can rely on the Order homotopy
theorem, Proposition~\ref{prop_segal}.
\end{pf}

As a consequence of Theorem~\ref{thm_teorema} we obtain the
following result answering the question from the title to this
section. Almost at the same time, actually a few days earlier,
this result was announced by P\' eter Csorba, \cite{Csorba2004}.

\begin{thm}{\rm (\cite{Csorba2004})}
For each finite, free $\mathbb{Z}_2$-complex $K$ there exists a
graph $G$ such that the associated graph complex $\Lambda(G)$ is
$\mathbb{Z}_2$-homotopy equivalent to $K$.
\end{thm}

\begin{pf}
The result is a consequence of Theorem~\ref{thm_teorema} since
obviously the face semilattice $\Phi(K)$ of $K$ is a free
$I$-semilattice with ``fat'' intervals.
\end{pf}

\bigskip

\section{Complexes $Hom(G,H)$ and the Lov\'{a}sz conjecture}\label{sec_homomorphisms}

\bigskip

The notion of a box poset associated to a $WI$-poset $(P,C)$,
Definition~\ref{def_box_poset}, admits several generalizations in
different directions. Here is one of the possibilities which
relates this construction to complexes $Hom(G,H)$. Recall that
these objects were introduced by L.~Lov\'{a}sz whose well known
conjecture about chromatic numbers of graphs $G$ with
$k$-connected complexes $Hom(C_{2r-1},G)$ was recently confirmed
by Bobson and Kozlov in \cite{BobsonKozlov2}.

\begin{defn}
\label{def_G_box_poset} Suppose that $G=(V_G,E_G)$ is a graph on
the ground set $[n]$, $V_G\subseteq [n]$. The $G$-\emph{box poset}
$G$-$\mathfrak{B}(P)$ associated to a $WI$-poset $(P,C)$ is a
subposet of $P^n$ defined by
\begin{equation*}
G{\text -}\mathfrak{B}(P) :=\left\{ (x_1,\ldots ,x_n)\in P^n \mid
(\forall i\neq j) \, \{i,j\}\in E_G \Rightarrow x_i\leq C(x_j)
\And x_j\leq C(x_i)\right\} .
\end{equation*}
\end{defn}
If $G=K_2$ is the complete graph on two vertices, the $G$-box
poset $G$-$\mathfrak{B}(P)$ reduces to the box poset
$\mathfrak{B}(P)$ from Section~\ref{sec_involutive}. More
importantly, if $P = P_H = \{B\subset V_H \mid
CN(B)\neq\emptyset\}$ is the $WI$-poset associated to a graph $H =
(V_H,E_H)$, then $G$-$\mathfrak{B}(P)$ is the face poset
associated to the polyhedral complex $Hom(G,H)$,
\cite{BobsonKozlov1} \cite{BobsonKozlov2}. The fact that
$Hom(K_2,G)$ is one of avatars of the graph complex $\Lambda(G)$,
Definition~\ref{def_lambda}, is already an indication of the
importance of the complex $Hom(G,H)$. Lov\'{a}sz conjectured that
if $Hom(C_{2r+1},G)$ is $k$-connected for some $r\geq 1$, where
$C_d$ is the $d$-cycle, then $\chi(G)\geq k+4$. This conjecture
was recently proved by Bobson and Kozlov, \cite{BobsonKozlov2}.
The proof is reasonably long and quite intricate involving a
variety of different techniques. In particular it required a
detailed combinatorial and homological analysis of polyhedral
complexes $Hom(G, K_n)$ with a special emphasis on the complex
$Hom(C_{2r+1},K_n)$, \cite{BobsonKozlov1} \cite{BobsonKozlov2}.
Having in mind that the existence of different models for the
graph complex $\Lambda(G)$ makes them more accessible, it is
interesting to ask if $Hom(G,H)$, and in particular the complex
$Hom(C_{2r+1},G)$, also have different incarnations. Even if the
answer is negative, it may be of some interest to establish a
``hierarchy theorem'' in the spirit of Theorem~1 in
\cite{MatZieg}.

Let us start with the observation that most of the complexes and
posets from
Sections~\ref{sec_all_equivalent}--\ref{sec_applications} do have
their analogs in the broader context of $Hom(G,H)$ complexes. We
will not attempt to give a complete analysis here. Instead, we
select some model cases and give examples which illuminate
potential use of these more general objects. For example the poset
$G$-$\mathfrak{B}_{ex}(P)$ is a relative of $G$-$\mathfrak{B}(P)$
obtained if in the Definition~\ref{def_G_box_poset} we allow some,
but not all entries in the vector $(x_1,\ldots,x_n)$ to be equal
to an added new minimum element $\hat{0}$. The following
definition is just a repetition of the definition of
$G$-$\mathfrak{B}_{ex}(P)$ in the case of the $WI$-poset $P=P_H$.
The notation emphasizes the fact that the new complexes are
relatives of the poset (complex) $Hom(G,H)$.

\begin{defn}
\label{def_homex} The extended $Hom$-poset $Hom_{ex}(G,H)$ is a
poset whose elements are all functions $\phi : V_G\rightarrow
2^{V_H}$, such that $\phi(i)\neq\emptyset$ for some $i\in V_G$,
for each edge $\{i,j\}\in E_G$, $\phi(i)\cap\phi(j)=\emptyset$ and
$\forall x\in\phi(i)\,\forall y\in\phi(j) \, \{x,y\}\in E_H$.
\end{defn}
The reader familiar with \cite{BobsonKozlov2} will notice right
away that our $Hom_{ex}(G,H)$ is nothing but the complex
$Hom_{+}(G,H)$ which plays a very important role in the analysis
leading eventually to the proof of Lov\'{a}sz conjecture! Note
that $Hom(G,H)$ is a subposet of $Hom_{ex}(G,H)$ and, as a
consequence of the analysis from \cite{BobsonKozlov2}, one cannot
expect that these two complexes are homotopy equivalent in
general. Note also that in the case of a complete graph $H=K_n$,
$\phi\in Hom_{ex}(G,K_n)$ iff $\phi(i)\neq\emptyset$ for some $i$
and $\phi(i)\cap\phi(j)=\emptyset$ for each edge $\{i,j\}\in E_G$.
In this case $Hom_{ex}(G,H)$ can be seen as a subposet
(subcomplex) of a join $(\Delta^{n-1})^{\ast V_G} =
\Delta^{n-1}\ast\ldots\ast\Delta^{n-1}$ of $\vert V_G\vert$-copies
of the $(n-1)$-simplex $\Delta^{n-1}$ spanned by vertices of the
graph $K_n$. In order to simplify notation, from here on we assume
that $V_G=[m], V_{K_n}=[n]$ and to each function $\phi\in
Hom_{ex}(G,K_n)$ we associate  its ``graph'' $\Gamma(\phi)\subset
[m]\times [n]$, where $\Gamma(\phi)\cap(\{i\}\times [n]) =
\{i\}\times\phi(i)$. In this notation, $Hom_{ex}(G,K_n)\subset
(\Delta^{n-1})^{\ast[m]}$. More importantly, the condition
$\phi(i)\cap\phi(j)=\emptyset$ for each edge $\{i,j\}\in E_G$
indicates that $Hom_{ex}(G,K_n)$ is, as a simplicial complex, a
$G$-deleted join of simplices $\Delta^{n-1}$ in the sense of the
following definition.

\begin{defn}
\label{def_G_del_join} Suppose that $G$ is a graph on $[n]$ as a
ground set, $V_G\subset [n]$. Let $\{K_i\}_{i=1}^n$ be a
collection of $n$-copies of a simplicial complex $K$. Then the
$G$-deleted join of $K$ is the simplicial subcomplex $K_G^{\ast
n}$ of $K\ast\ldots\ast K = K^{\ast n}$, where
$\theta_1\ast\ldots\ast\theta_n\in K_G^{\ast n}$ iff
$\theta_i\cap\theta_j=\emptyset$ for each edge $\{i,j\}\in E_G$.
\end{defn}
If $G=K_m$ is a complete graph then $G$-deleted join $K_G^{\ast m}
= K_{\Delta}^{\ast m}$ reduces to the usual deleted join operation
of simplicial complexes, \cite{Matousek} Section~5.5,
\cite{Sarkaria}, \cite{guides}. The well known relation $(K\ast
L)_{\Delta}^{\ast k} \cong K_{\Delta}^{\ast k}\ast
L_{\Delta}^{\ast k}$ easily generalizes to the following result
\begin{lemma}
\label{lema_well_known} Suppose that $K$ and $L$ are simplicial
complexes and let $G=(V_G,E_G)$ be a graph on the ground set
$[m]$, $V_G\subset [m]$. Then,
\[
(K\ast L)_{G}^{\ast m} \cong K_{G}^{\ast m}\ast L_{G}^{\ast m}.
\]
\end{lemma}
An immediate consequence of Lemma~\ref{lema_well_known} is the
relation
\begin{equation}\label{eqn_indep}
Hom_{ex}(G,K_n)= (\Delta^{n-1})_G^{\ast m}\cong ((pt)^{\ast
n})_G^{\ast m} \cong ((pt)_G^{\ast m})^{\ast n} .
\end{equation}
The complex $(pt)_G^{\ast m}$ is well known as the complex
$Ind(G)$ of all independent sets in a graph $G$. Hence the
equation (\ref{eqn_indep}) is nothing but (a half of) the
Proposition~3.2. from \cite{BobsonKozlov2} in disguise. This shows
that the study of complexes $Hom_{ex}(C_m,K_n)$ is reduced to the
study of complexes $Ind(C_n)$, which is the subject of our next
section.

\bigskip

\section{$Ind(L_n)$ and $Ind(C_n)$ as $\mathbb{Z}_2$-complexes}
\label{sec_ind}

\bigskip

Let us denote by $\one$ and $\varepsilon$ respectively the trivial
and nontrivial real representations of $\mathbb{Z}_2$. Given an
Euclidean vector space $V$, let $S(V)$ be the associated unit
sphere. If $V$ is an orthogonal representation of $\mathbb{Z}_2$,
the sphere $S(V)$ is a $\mathbb{Z}_2$-space. For example $S(\one)$
and $S(\varepsilon)$ are both $2$-element sets, the first with
trivial and the second with non-trivial action of $\mathbb{Z}_2$.
Recall the well known fact that $S(U\oplus V)\cong S(U)\ast S(V)$.
For example if $V=p\one\oplus q\varepsilon$ then $S(V)$ is the
sphere in $\mathbb{R}^{p+q}$ equipped with the action of
$\mathbb{Z}_2 = \{1,\omega\}$ such that $\omega(x_1,\ldots
,x_p,y_1,\ldots ,y_q) = \omega(x_1,\ldots ,x_p,-y_1,\ldots
,-y_q)$.

\begin{defn}
Define graphs $L_n$ and $C_n$ on $[n]=\{1,\ldots, n\}$ as the
ground set by the conditions
\[
(i,j)\in E_{L_n} \Leftrightarrow \vert i - j\vert =1 \enskip
\text{ {\rm and} } \enskip (i,j)\in E_{C_n} \Leftrightarrow \vert
i - j\vert =1 \text{ \rm{({\rm mod} $n$)}}.
\]
Given an interval $[p,q]$ in $[n]$, let $L_{[p,q]}\cong L_{q-p+1}$
be the complete subgraph of $L_n$ on $[p,q]$ as the set of
vertices. Define $\mathbb{Z}_2$-actions on both $L_n$ and $C_n$ by
the involution $\omega : [n]\rightarrow [n]$ which sends $i$ to
$n-i+1$. Let $Ind(L_n)$ and $Ind(C_n)$ be the associated complexes
of independent sets with inherited $\mathbb{Z}_2$-actions.
\end{defn}

Homotopy types of spaces $Ind(L_n)$ and $Ind(C_n)$ were determined
in \cite{Kozlov99}. The question of finding the associated
$\mathbb{Z}_2$-homotopy types appeared as a natural step in the
approach of Bobson and Kozlov to the solution of Lov\'{a}sz
conjecture, notably in the evaluation of the height of the first
Stiefel-Whitney class of the $\mathbb{Z}_2$-complex $Hom(C_{2r+1},
K_n)$, \cite{BobsonKozlov2} Sections~2.2 and 4.1. Their methods
permitted them to evaluate only the homotopy types of the
associated orbit spaces $Ind(C_n)/\mathbb{Z}_2$ but this turned
out to be sufficient for the intended application.

In this section we strengthen this result of Bobson and Kozlov by
demonstrating how the $\mathbb{Z}_2$-homotopy types of these
complexes can be determined, again relying on the Bredon's
theorem. We restrict ourselves to the analysis of two important
special cases. In the other cases, corresponding to other values
of $n$, the proofs are similar in spirit and rely on similar
ideas.

\begin{prop}
\label{prop_Ind1}{\rm (\cite{BobsonKozlov2})} Assume that $n=6p-1$
and let $S(\one)$ and $S(\varepsilon)$ be $2$-element sets {\rm
(}$0$-dimensional spheres{\rm )} respectively with trivial and
non-trivial action of $\mathbb{Z}_2$. Then $Ind(L_n)$ is a
$\mathbb{Z}_2$-complex which is $\mathbb{Z}_2$-homotopy equivalent
to the join of $p$ copies of $S(\one)$ and $p$ copies of
$S(\varepsilon)$,
\begin{equation}
\label{eqn:p_and_p} Ind(L_{6p-1}) \simeq_{\mathbb{Z}_2}
S(\one)^{\ast p}\ast S(\varepsilon)^{\ast p}.
\end{equation}
\end{prop}

Before we commence the proof of the proposition let us introduce
some auxiliary definitions and useful lemmas. Given a graph $G =
(V_G, E_G)$ and a subset $K\subset V_G$, define $G\setminus K$ as
the graph obtained from $G$ by removing $K$ and all edges incident
to vertices in $K$. For example $G\setminus v$ is obtained from
$G$ by removing a vertex $v$, while $G\setminus St(v)$ is the
graph obtained from $G$ if $K = St(v)$ is the star of $v$, $St(v)
= \{v\}\cup \{w\in V_G \mid (v,w)\in E_G\}$.

\begin{lemma}
\label{lema_decomp} There is a decomposition $Ind(G) = X\cup Y$
where $X = Ind(G\setminus v)$ and $Y = \{v\}\ast Ind(G\setminus
St(v))$ where $X\cap Y = Ind(G\setminus St(v))$.
\end{lemma}

\begin{lemma}
\label{lema_basic} Suppose that $a,b,v\in V_G$ are three distinct
vertices in a graph $G=(V_G,E_G)$ such that both $(a,b)\in E_G$
and $(b,v)\in E_G$. Moreover we assume that $a$ is not connected
with any other vertex in $G$, i.e. ${\rm deg}(a)=1$. Then
$Ind(G)\simeq Ind(G\setminus v)$.
\end{lemma}

\begin{pf}
Note that $Y$ in the decomposition $Ind(G) = X\cup Y$ in
Lemma~\ref{lema_decomp} is contractible, being a cone with vertex
$v$. The space $X\cap Y$ is also a cone since by assumption
$G\setminus St(v)$ has an isolated vertex $a$. The proof is
completed by invoking an easily established fact that if both $Y$
and $X\cap Y$ are contractible complexes then $X\cup Y\simeq X$.
\end{pf}

\medskip\noindent
{\bf Proof of Proposition~\ref{prop_Ind1}.} By successive
applications of Lemma~\ref{lema_decomp}, we are able to remove all
vertices from the set $K = \{3,6,\ldots ,3p,\ldots, 6p-3\}$
without changing the homotopy type of $Ind(L_{6p-1})$. In other
words, $Ind(L_{6p-1})\simeq Ind(L_{6p-1}\setminus K)$. Let us show
that the inclusion map $e : Ind(L_{6p-1}\setminus K)\rightarrow
Ind(L_{6p-1})$ is actually a $\mathbb{Z}_2$-homotopy equivalence.
By Theorem~\ref{thm_bredon}, it is sufficient to show that the
inclusion map $e^{\mathbb{Z}_2} : Ind(L_{6p-1}\setminus
K)^{\mathbb{Z}_2}\rightarrow Ind(L_{6p-1})^{\mathbb{Z}_2}$ of the
associated spaces of fixed points is also a homotopy equivalence.
Both $Ind(L_{6p-1})$ and $Ind(L_{6p-1}\setminus K)$ are
subcomplexes of the simplex $\Sigma$ spanned by vertices
$1,\ldots, n$. Identifying $\Sigma$ with its geometric realization
$\vert\Sigma\vert$, assume that vertices of $\Sigma$ are points
$v_1,\ldots, v_n$ in some vector space $V$. Note that $\Sigma$ is
also a $\mathbb{Z}_2$-space with the linear action which is on
vertices defined by $\omega(v_i) = v_{n+1-i}$. It is not difficult
to check that
\begin{equation}
x\in \Sigma^{\mathbb{Z}_2} \enskip \Leftrightarrow \enskip x =
t_1\frac{v_1+v_{6p-1}}{2} +\ldots +
t_{3p-1}\frac{v_{3p-1}+v_{3p+1}}{2} + t_{3p}x_{3p}
\end{equation}
where $t_j\geq 0$ and $\Sigma t_j=1$. We conclude that
$\Sigma^{\mathbb{Z}_2}$ is a simplex isomorphic to the face
$\Sigma_1$ of $\Sigma$ spanned by the vertices $\{v_1,v_2,\ldots,
v_{3p}\}$, where the isomorphism $I : \Sigma_1\rightarrow\Sigma$
is the linear extension of the map $v_i\mapsto(v_i+v_{6p-i})/2$.
The fixed point spaces $Ind(L_{6p-1}\setminus K)^{\mathbb{Z}_2}$
and $Ind(L_{6p-1})^{\mathbb{Z}_2}$ are subspaces of
$\Sigma^{\mathbb{Z}_2}$ which can be viewed, via isomorphism $I$,
as subspaces of $\Sigma_1$. It immediately follows that
$Ind(L_{6p-1})^{\mathbb{Z}_2}\cong Ind(L_{[1,3p]})\cong
Ind(L_{3p})$ and $Ind(L_{6p-1}\setminus K)^{\mathbb{Z}_2}\cong
Ind(L_{[1,3p]}\setminus K')$ where $K'=\{3,6,\ldots, 3p\}$. Again,
by applications of Lemma~\ref{lema_basic} and successive removal
of vertices in $K'$, we conclude that the inclusion map
$Ind(L_{[1,3p]}\setminus K')\hookrightarrow Ind(L_{[1,3p]})$ is a
homotopy equivalence, hence $e^{\mathbb{Z}_2} :
Ind(L_{6p-1}\setminus K)^{\mathbb{Z}_2}\rightarrow
Ind(L_{6p-1})^{\mathbb{Z}_2}$ is a homotopy equivalence. This, in
light of Theorem~\ref{thm_bredon}, completes the proof that $e$ is
a $\mathbb{Z}_2$-homotopy equivalence. Note that
$Ind(L_{6p-1}\setminus K)$ is isomorphic to the following join of
circles,
\[
Ind(L_{6p-1}\setminus K)\cong Ind(L_{[1,2]}\cup
L_{[6p-2,6p-1]})\ast\ldots\ast Ind(L_{[3p-2,3p-1]}\cup
L_{[3p+1,3p+2]}).
\]
This, together with the fact that
\[
Ind(L_{[j,j+1]}\cup L_{[6p-j-1,6p-j]})\cong S(\one)\ast
S(\varepsilon)
\]
finally completes the proof of Proposition~\ref{prop_Ind1}. \quad
$\blacksquare$

\bigskip
Let $K$ be a finite simplicial complex and assume that $\sigma\in
K$ is a simplex which is maximal in the sense that it is not a
proper face of any other simplex $\tau\in K$. Let ${\rm
dim}(\sigma)=k$. If $K':=K\setminus \{\sigma\}$ then the geometric
realization of $K'$ is obtained from the geometric realization of
$K$ by removing the interior $\stackrel{\circ}{\sigma}$ of
$\sigma$, $\vert K'\vert = \vert K\vert\setminus
\stackrel{\circ}{\sigma}$. If $K'$ is contractible then $K\simeq
K/K'\cong\sigma/\partial\sigma\cong S^k$. In this case we call
$\sigma$ a {\em generating simplex} of $K$. Of course, it is not
true that a complex homotopy equivalent to a sphere must have a
generating simplex. For example $S^{k-1}\times I$ is a pure
$k$-dimensional complex homotopic to $S^{k-1}$ which consequently
cannot have a ($(k-1)$-dimensional) generating simplex. The
following lemma gives a sufficient condition for the existence of
generating simplices.

\begin{lemma}
\label{lema_gen} Let $K$ be a finite simplicial complex and assume
that $L\subset K$ is a subcomplex of $K$ simplicially isomorphic
to a triangulation of a $k$-sphere $S^k$. Assume that the
inclusion map $e: L\rightarrow K$ is a homotopy equivalence and
let $\sigma\in L$ be a $k$-simplex which is maximal in $K$, i.e.\
such that $\sigma$ is not a proper face of a simplex $\tau\in K$.
Then $\sigma$ is a generating simplex for $K$ in the sense that
the complex $K'=K\setminus\{\sigma\}$ is contractible and
$K/K'\simeq \sigma/\partial\sigma\cong S^k$.
\end{lemma}
\begin{pf}
By assumption $L$ is a weak deformation retract of $K$ hence a
strong deformation retract, \cite{Spanier} Section~I.4. Since
$\sigma$ is maximal in $K$ we observe that $L\setminus\{\sigma\}$
is a strong deformation retract of $K\setminus\{\sigma\}$. Since
$L\setminus\{\sigma\}$ is contractible, $K\setminus\{\sigma\}$ is
also contractible and the result follows.
\end{pf}

\begin{exmp}
\label{ex_gen} The proof of Proposition~\ref{prop_Ind1} reveals
that the complexes $Ind(L_{6p-1})$ and $L = Ind(L_{6p-1}\setminus
K)$ satisfy the conditions of Lemma~\ref{lema_gen}. A simplex
$\sigma$ in $Ind(L_{6p-1}\setminus K)$ is maximal in
$Ind(L_{6p-1})$ if and only if there are at most two vertices from
the ground set $[6p-1]$ separating two consecutive vertices in
$\sigma$. Hence an example of a generating simplex is
 \[
 \tau = \{2,5,\ldots ,3p-4, 3p-1,3p+1,3p+4,\ldots, 6p-5,6p-2\}.
 \]
 \end{exmp}

\begin{prop}
\label{prop_Ind2} The complex $Ind(C_{6p-1})$ is
$\mathbb{Z}_2$-homotopy equivalent to the sphere $S^{2p-1}\subset
\mathbb{R}^{2p}$ with the action of $\mathbb{Z}_2=\{1,\omega\}$
given by the formula
\[\omega(x_1,\ldots, x_p,y_1,\ldots, y_p) = (x_1,\ldots,
x_p,-y_1,\ldots, -y_p).\] In other words,
\begin{equation}
\label{eqn_ind2} Ind(C_{6p-1})\simeq_{\mathbb{Z}_2} S(\one)^{\ast
p}\ast S(\varepsilon)^{\ast p}.
\end{equation}
\end{prop}
\begin{pf}
It was show in \cite{Kozlov99} that $Ind(C_{6p-1})\simeq
S^{2p-1}$. Moreover it was shown that the simplex
\[
 \sigma = \{2,5,\ldots ,3p-4, 3p-1,3p+1,3p+4,\ldots, 6p-5,6p-2\}
\]
is a generating simplex for the complex $Ind(C_{6p-1})$. The
reader is invited to prove this fact along the lines of proofs of
Proposition~\ref{prop_Ind1} and Example~\ref{ex_gen}. As a
consequence we know that $Int(C_{6p-1})\setminus
\stackrel{\circ}{\sigma}$ is contractible. Let us note that
$\sigma$ is $\mathbb{Z}_2$-invariant with respect to the
$\mathbb{Z}_2$-action on $Ind(C_{6p-1})$ which, as we recall,
arises from the involution $\omega : [n]\rightarrow [n],\,
\omega(j):= n+1-j$. It follows that $\sigma/\partial\sigma\cong
S^k$ is a $\mathbb{Z}_2$-space and there is an obvious
$\mathbb{Z}_2$-equivariant collapsing map $f :
Ind(C_{6p-1})\rightarrow \sigma/\partial\sigma$. Let us show that
this map is a $\mathbb{Z}_2$-homotopy equivalence. Since $f$ is a
homotopy equivalence, by Theorem~\ref{thm_bredon} it is sufficient
to show that $f^{\mathbb{Z}_2} :
Ind(C_{6p-1})^{\mathbb{Z}_2}\rightarrow
(\sigma/\partial\sigma)^{\mathbb{Z}_2}$ is a homotopy equivalence.
We follow the same strategy as in the proof of
Proposition~\ref{prop_Ind1}, in particular we use the map $I$ to
relate the fixed point sets to (subspaces) of independence
complexes. For example, as in the proof of
Proposition~\ref{prop_Ind1}, $Ind(C_{6p-1})^{\mathbb{Z}_2}$ is
isomorphic to the complex $Ind(L_{[1,3p]})$. Similarly,
$(\sigma/\partial\sigma)^{\mathbb{Z}_2}$ is isomorphic to the
space $\sigma_0/\partial\sigma_0$, where $\sigma_0$ is the simplex
in $Ind(L_{[1,3p]})$ spanned by vertices $\{2,5,\ldots,
3p-4,3p-1\}$. Note that $\sigma_0$ is a maximal simplex in
$Ind(L_{[1,3p]})$. Moreover $\sigma_0$ is a $(p-1)$-dimensional
simplex in the $(p-1)$-sphere
\[
S^{p-1}\cong Ind(L_{[1,2]})\ast\ldots\ast Ind(L_{[3p-2,3p-1]})
\]
which is a deformation retract of $Ind(L_{[1,3p]})$. Hence
$\sigma_0$ is a generating simplex in the complex
$Ind(L_{[1,3p]})$ which shows that the collapsing map $f_0:
Ind(L_{[1,3p]})\rightarrow\sigma_0/\partial\sigma_0$ is a homotopy
equivalence. This in turn implies that $f^{\mathbb{Z}_2} :
Ind(C_{6p-1})^{\mathbb{Z}_2}\rightarrow
(\sigma/\partial\sigma)^{\mathbb{Z}_2}$ is a homotopy equivalence
and by Bredon's theorem $Ind(C_{6p-1})$ is $\mathbb{Z}_2$-homotopy
equivalent to the $\mathbb{Z}_2$-space $\sigma/\partial\sigma$.

In order to determine the $\mathbb{Z}_2$-structure of the
$\mathbb{Z}_2$-space $\sigma/\partial\sigma$, note that it was
already done in the proof of Proposition~\ref{prop_Ind1}. Indeed,
the simplex $\sigma$ was shown there to be a generating simplex of
the complex $Ind(L_{6p-1}\setminus K)$ and the collapsing map
$Ind(L_{6p-1}\setminus K)\rightarrow \sigma/\partial\sigma$ is a
$\mathbb{Z}_2$-homotopy equivalence, again by an applications of
Bredon's theorem. This finally establishes the decomposition
(\ref{eqn_ind2}).
\end{pf}

\begin{rem}
The fact that both $Ind(L_{6p-1})$ and $Ind(C_{6p-1})$ have
identical $\mathbb{Z}_2$ decompositions,
Propositions~\ref{prop_Ind1} and \ref{prop_Ind2}, is not an
accident. Given a graph $G=(V_G,E_G)$ and an edge $e=(u,v)\in
E_G$, define $G\setminus e$ and $G\setminus St(e)$ as the graphs
$G\setminus e = (V_G, E_G\setminus\{e\})$ and $G\setminus St(e) =
G\setminus\{u,v\}$. Then there is a decomposition, cf.\ {\rm
\cite{Meshulam}}, $Ind(G\setminus e)=Ind(G)\cup \{u,v\}\ast
Ind(G\setminus St(e))$ where $Ind(G)\cap \{u,v\}\ast
Ind(G\setminus St(e))\cong Ind(G\setminus St(e))$. If \
$G=C_{6p-1}$ and $e=(1,6p-1)$, then
\[
Ind(G\setminus St(e)) =
Ind([6p-1]\setminus\{1,2,6p-2,6p-1\})\cong Ind(L_{6p-5})
\]
is contractible. It follows, along the lines of the proof of
Lemma~\ref{lema_basic}, that the natural inclusion map $e:
Ind(C_{6p-1})\rightarrow Ind(L_{6p-1})$ is a homotopy equivalence.
The map $e$ is $\mathbb{Z}_2$-equivariant and a repetition of the
argument already used in the proofs of
Propositions~\ref{prop_Ind1} and \ref{prop_Ind2} allows us to
conclude that $e$ is a $\mathbb{Z}_2$-homotopy equivalence.
\end{rem}

\bigskip

\bigskip

\end{document}